\documentclass{article}

\usepackage{amsmath, amsfonts, enumerate, amsthm}

\title{Grothendieck ring class of Banana and Flower graphs}
\author{Pedro Morales-Almaz\'an\\ \texttt{Pedro\_Morales@baylor.edu}}
\date{}

\newtheorem{thm}{Theorem}[section]
\newtheorem{pro}{Proposition}[section]
\newtheorem{cor}[thm]{Corollary}
\newtheorem{lem}[thm]{Lemma}

\begin{document}
\maketitle
\begin{abstract}
We define a special type of hypersurface varieties inside $\mathbb{P}_k^{n-1}$ arising from connected planar graphs and then find their equivalence classes inside the Grothendieck ring of projective varieties. Then we find a characterization for graphs in order to define irreducible hypersurfaces in general.
\end{abstract}

\section{Introduction}

This is result from a short communication session from the summer school Applications of Algebra and Topology in Quantum Field Theory held in Villa de Leyva, Colombia during the summer of 2011.

The first three sections are mainly based on \cite{a08}. In the first section we introduce the concept of the Grothendieck ring of varieties together with some important results in this ring such as the inclusion-exclusion principle. In the second section we define a graph polynomial and a graph hypersurface and find some special properties about such polynomials.
The third section is mainly concerned with the Grothendieck class of some special type of graphs, namely star, flower, polygon and banana graphs. In section four we explore graph hypersurfaces for graphs in general and find a necessary and sufficient condition for a graph to produce an irreducible graph hypersurface.

\section{Grothendieck Ring of Varieties}

The Grothendieck ring of varieties can be thought as a generalization of the Euler characteristic, \cite{s08}. It can be thought of as the quotient of the free abelian group generated by isomorphism classes $[X]$ of $k$-varieties, modulo inclusion-exclusion type relations, together with the product of two classes given by the class of the product, $\cite{a08, s08}$.

We start off with a characteristic zero field $k$ and the the category of algebraic varieties $\mathcal{V}_k$ defined over it. Then we look at the abelian group $K_0(\mathcal{V}_k)$ of isomorphism classes $[X]$ of varieties $X$ over the field $k$ with the relations
\begin{equation}\label{quo}
[X]-[Y]=[X\setminus Y],
\end{equation}
where $Y\subset X$ is a closed subvariety of $X$. Notice that setting $X=Y$ gives that the Grothendieck class of the empty set is zero,
\begin{equation}\label{zeroclass}
[\emptyset]=0\,.
\end{equation}
$K_0(\mathcal{V}_k)$ is called the \emph{Grothendieck group of varieties} which then can be turned into a ring by defining the product of two isomorphism classes by
\begin{equation}\label{pro}
[X]\cdot[Y]=[X\times Y].
\end{equation}

This can be also thought of as the quotient of the free abelian group generated by the symbols $[X]$ by the relation (\ref{quo}) and with the product (\ref{pro}), \cite{s08}.

The use of the Grothendieck ring of varieties is very useful when looking for invariants. We call an \emph{additive invariant} a map $\chi:\mathcal{V}_k\to R$, with values on a commutative ring $R$, that satisfies
\begin{enumerate}[1.]
\item $\chi(X)=\chi(Y)$ if $X$ and $Y$ are isomorphic
\item $\chi(X\setminus Y)=\chi(X)-\chi(Y)$, for $Y\subset X$ closed
\item $\chi(X\times Y)=\chi(X)\chi(Y)$
\end{enumerate}

Thus, an additive invariant is the same as a ring homomorphism $\chi$ from the Grothendieck ring of varieties to the ring R, \cite{a08}. Some interesting additive invariants are the Euler characteristic, as mentioned before, and the Hodge polynomial, \cite{a08, s08}.

For the following discussion, let $\mathbb{L}=[\mathbb{A}_k^1]$ be the equivalence class in the Grothendieck ring of varieties of the one dimensional affine space $\mathbb{A}_k^1$. From the product (\ref{pro}) of $K_0(\mathcal{V}_k)$ we have that the Grothendieck class of the affine spaces $\mathbb{A}_k^n$ are given in terms of $\mathbb{L}$ as
\begin{equation}
[\mathbb{A}_k^n]=\mathbb{L}^n\,.
\end{equation}

Also, let $\mathbb{T}$ be the class of the multiplicative group $k^\times$, which is just $\mathbb{A}_k^1$ without a point. Thus we have that
\begin{equation}\label{tor}
\mathbb{T}=[\mathbb{A}_k^1]-[\mathbb{A}_k^0]=\mathbb{L}-1,
\end{equation}
where $1$ is the class of a point in $K_0(\mathcal{V}_k)$.

$\mathbb{L}$ and $\mathbb{T}$ are useful to find the Grothendieck class of a variety together with the inclusion-exclusion principle

\begin{thm}[Inclusion-Exclusion Principle]

Let $X, Y$ be varieties over $k$, then
\begin{equation}
[X\cup Y]=[X]+[Y]-[X\cap Y]
\end{equation}

\begin{proof}

Since we have that $X\setminus (X\cap Y)=(X\cup Y)\setminus Y$, then
\begin{equation}
[X]-[X\cap Y]=[X\setminus (X\cap Y)]=[(X\cup Y)\setminus Y]=[X\cup Y]-[Y]\,,
\end{equation}
and hence the result follows.
\end{proof}
\end{thm}

From this we find as a corollary the classes the projective spaces,

\begin{cor}
We have that
\begin{equation}\label{projclass}
 [\mathbb{P}_k^n]=\frac{\mathbb{L}^{n+1}-1}{\mathbb{L}-1}\,,
\end{equation}
where the fraction is taken as a short hand notation for the corresponding summation.

\begin{proof}
Since $\mathbb{P}_k^{n+1}\setminus \mathbb{A}_k^{n+1}\simeq \mathbb{P}_k^n$ for all $n\in\mathbb{N}$, by induction and the previous result, we have that
\begin{equation}
[\mathbb{P}_k^n]=1+\mathbb{L}+\mathbb{L}^2+\dots+\mathbb{L}^n\,,
\end{equation}
from which the result follows.
\end{proof}
\end{cor}

\section{Graphs}

\subsection{Graph Hypersufaces}
In this section, let $\Gamma$ be a connected planar graph with $n$ edges and label them with the variables $t_1,t_2,\dots, t_n$ of the polynomial ring $Z[t_1,t_2,\dots, t_n]$. With this, define the graph polynomial associated to $\Gamma$ by
\begin{equation}
\Psi_\Gamma(t)=\sum_{T\subset\Gamma}\prod_{e\notin E(T)} t_e\,,
\end{equation}
where $T$ is runs through all spanning trees of $\Gamma$ and $E(T)$ is the set of all edges of $T$.

For a graph with $v$ vertices, the pigeon hole principle gives that all spanning trees have $v-1$ edges, as otherwise there would be a cycle. Hence the polynomial $\Psi_\Gamma(t)$ is a homogeneous polynomial of degree $n-v+1$ and one can define the graph hypersurface associated to the graph $\Gamma$ by
\begin{equation}
X_\Gamma=\{t=(t_1:t_2:\dots:t_n)\in\mathbb{P}_k^{n-1}|\Psi_\Gamma(t)=0\}\,.
\end{equation}
Since $\Psi_\Gamma(t)$ is homogeneous, $X_\Gamma$ is well defined as a projective variety.

\subsection{Dual Graph}
Given a planar connected graph $\Gamma$, define the dual graph $\Gamma^\vee$ by the following:
\begin{enumerate}[1.]
\item embed $\Gamma$ in $\mathbb{S}^2$
\item for each region defined by $\Gamma$ on $\mathbb{S}^2$ assign a vertex of $\Gamma^\vee$
\item if two regions share an edge, join the corresponding vertices on $\Gamma^\vee$
\end{enumerate}

This definition depends on the particular embedding used. Different embeddings of the same graph $\Gamma$ might lead to different dual graphs $\Gamma^\vee$, but the resulting graph polynomials $\Psi_{\Gamma^\vee}(t)$ are the same up to relabeling. This fact follows from the relation between the graph polynomial of $\Gamma$ and the graph polynomial of the dual $\Gamma^\vee$ via the Cremona transformation described below, \cite{a08,b07}. This makes the graph hypersurface of the dual graph $X_{\Gamma^\vee}$ a well defined object.

\subsection{Cremona Transformation}

The Cremona transformation on $\mathbb{P}_k^{n-1}$ is given by
\begin{gather}
\mathcal{C}:\mathbb{P}_k^{n-1}\to\mathbb{P}_k^{n-1}\\
(t_1:t_2:\dots:t_n)\mapsto\left(\frac{1}{t_1}:\frac{1}{t_2}:\dots:\frac{1}{t_n}\right)
\end{gather}
which is well defined outside the coordinate hyperplanes
\begin{equation}\label{sigman}
\Sigma_n=\left\{(t_1:t_2:\dots:t_n)\in\mathbb{P}_k^{n-1}|\prod_{i=1}^nt_i=0\right\}\,.
\end{equation}
This Cremona transform is useful to relate the graph hypersurfaces of a graph $\Gamma$ and its dual $\Gamma^\vee$. For this, consider the graph of the Cremona transform in $\mathbb{P}_k^{n-1}\times\mathbb{P}_k^{n-1}$ and define $\mathcal{G}(\mathcal{C})$ to be its closure. Then we have the two projections $\pi_1$ and $\pi_2$ of $\mathcal{G}(\mathcal{C})$ into $\mathbb{P}_k^{n-1}\times\mathbb{P}_k^{n-1}$ related by
\begin{equation}\label{projections}
\mathcal{C}\circ\pi_1=\pi_2\,.
\end{equation}

$\mathcal{G}(\mathcal{C})$ is then a closed subvariety of $\mathbb{P}_k^{n-1}\times\mathbb{P}_k^{n-1}$ with equations (see \cite{a08}, Lemma 1.1)
\begin{equation}
t_1t_{n+1}=t_2t_{n+2}=\dots=t_nt_{2n}\,,
\end{equation}
where we used the coordinates $(t_1:t_2:\dots:t_{2n})$ for $\mathbb{P}_k^{n-1}\times\mathbb{P}_k^{n-1}$.

Now, we have a relation between the graph polynomials of a graph $\Gamma$ and its dual $\Gamma^\vee$ given by

\begin{pro}\label{dualpol}
Let $\Gamma$ be a graph with $n$ edges. Then the graph polynomials of $\Gamma$ and its dual $\Gamma^\vee$ are related by
\begin{equation}
\Psi_\Gamma(t)=\left(\prod_{i=1}^nt_i\right)\Psi_{\Gamma^\vee}\left(\frac{1}{t}\right)\,,
\end{equation}
and hence the corresponding graph hypersurfaces of $\Gamma$ and $\Gamma^\vee$ are related via the Cremona transformation by
\begin{equation}
\mathcal{C}\left(X_{\Gamma}\cap\left(\mathbb{P}_k^{n-1}\setminus\Sigma_n\right)\right)=X_{\Gamma^\vee}\cap\left(\mathbb{P}_k^{n-1}\setminus\Sigma_n\right)
\end{equation}

\begin{proof}
The proof follows from the definition of the graph polynomial and the combinatorial properties between $\Gamma$ and $\Gamma^\vee$ (see \cite{a08}, Lemma 1.3).
\end{proof}
\end{pro}

This result also gives an isomorphism between the graph hypersurfaces of $\Gamma$ and $\Gamma^\vee$ away from the coordinate hyperplanes which can be summarized in the following

\begin{cor}\label{coro}
Given a graph $\Gamma$, we have that the graph hypersurface of its dual is
\begin{equation}
X_{\Gamma^\vee}=\pi_2(\pi_1^{-1}(X_\Gamma))\,,
\end{equation}
where $\pi_i$ are the projections given on (\ref{projections}). Also, we have that the Cremona transform gives a (biregular) isomorphism
\begin{equation}
\mathcal{C}:X_\Gamma\setminus\Sigma_n\to X_{\Gamma^\vee}\setminus \Sigma_n\,,
\end{equation}
and the projection $\pi_2:\mathcal{G}(\mathcal{C})\to\mathbb{P}_k^{n-1}$ restricts to an isomorphism
\begin{equation}
\pi_2:\pi_1^{-1}\left(X_\Gamma\setminus\Sigma_n\right)\to X_{\Gamma^\vee}\setminus\Sigma_n\,.
\end{equation}

\begin{proof}
The result follows from Proposition \ref{dualpol} and the geometric properties of $\mathcal{G}(C)$ (see \cite{a08}, Lemma 1.2, Lemma 1.3, Corollary 1.4).
\end{proof}
\end{cor}

\section{Banana and Flower graph hypersurfaces}

\subsection{Flower graphs}
Let us start with the simplest of the graphs that define a graph hypersurface. Let $\Gamma$ be a star graph with $n$ edges. Then \
\begin{equation}
\Psi_\Gamma(t)=1\,,
\end{equation}
and hence
\begin{equation}
X_\Gamma=\emptyset\,.
\end{equation}
On the other hand, the dual graph $\Gamma^\vee$ is a flower graph consisting with only one vertex and $n$ loops. Therefore the graph polynomial associated to $\Gamma^\vee$ is given by
\begin{equation}
\Psi_{\Gamma^\vee}(t)=t_1t_2\cdots t_n\,,
\end{equation}
which shows that the graph hypersurface is
\begin{equation}
X_{\Gamma^\vee}=\{t\in\mathbb{P}_k^{n-1}|\Psi_{\Gamma^\vee}(t)=0\}=\Sigma_n\,.
\end{equation}
Therefore, by means of Corollary \ref{coro}, we have an isomorphism between the graph hypersurfaces of the star graph and the flower graph away from the coordinate hyperplanes, which is indeed the case,
\begin{equation}
X_\Gamma\setminus \Sigma_n=\emptyset=X_{\Gamma^\vee}\setminus \Sigma_n\,.
\end{equation}

From \ref{zeroclass} we have that the Grothendieck class of $X_\Gamma$ is
\begin{equation}
[X_\Gamma]=0\,.
\end{equation}
Now, for $X_{\Gamma^\vee}$, recall that from \ref{sigman}, we have that the complement of $\Sigma_n$  consists of all tuples where all the elements are different from zero, which is a copy of $(k^\times)^{n-1}$ and gives the class of $\mathbb{T}^{n-1}$.
\begin{pro}
The Grothendieck class of $\Sigma_n$ is given by
\begin{equation}
[X_{\Gamma^\vee}]=[\Sigma_n]=\frac{(1+\mathbb{T})^n-1-\mathbb{T}^n}{\mathbb{T}}=\sum_{i=1}^{n-1}\binom{n}{i}\mathbb{T}^{n-1-i}\,.
\end{equation}
\begin{proof}
Since $X_{\Gamma^\vee}=\Sigma_n=\mathbb{P}_k^{n-1}\setminus (k^\times)^{n-1}$, the result follows from the definition of Grothendieck class, the relations given in (\ref{projclass}) and (\ref{tor}).
\end{proof}
\end{pro}

\subsection{Banana graphs}

The next interesting type of graphs that give rise to graph hypersurfaces are the polygons. Let $\Gamma$ be a polygon  with $n$ edges, i.e. a graph with $n$ edges, $n$ vertices in which each vertex has degree 2. Here, the graph polynomial is then given by
\begin{equation}
\Psi_\Gamma(t)=t_1+t_2+\dots+t_n\,,
\end{equation}
and hence the graph hypersurface is the hyperplane
\begin{equation}
X_\Gamma=\{t\in\mathbb{P}_k^{n-1}|t_1+t_2+\dots+t_n=0\}=:\mathcal{L}\,.
\end{equation}
Notice that choosing any $n-1$ points for $t_1,t_2,\dots,t_{n-1}$ in $\mathcal{L}$ gives a unique value for $t_n$, and hence we can identify $\mathcal{L}$ with $\mathbb{P}_k^{n-2}$. Thus we find the Groethendieck class of $\mathcal{L}$ to be given by
\begin{equation}
[\mathcal{L}]=[X_{\Gamma}]=[\mathbb{P}_k^{n-2}]=\frac{\mathbb{L}^{n-1}-1}{\mathbb{L}-1}=\frac{(\mathbb{T}+1)^{n-1}-1}{\mathbb{T}}\,.
\end{equation}

The dual graph for $\Gamma$ is a banana graph, which consists of 2 vertices of degree $n$ each, and $n$ edges joining them. The corresponding graph polynomial is then given by the $(n-1)-$th elementary symmetrical polynomial in $t_i$,
\begin{equation}
\Psi_{\Gamma^\vee}(t)=t_2t_3\cdots t_{n}+t_1t_3\cdots t_n+\dots+t_1t_2\cdots\hat{t_i}\cdots t_n+\dots+t_1t_2\cdots t_{n-1}\,,
\end{equation}
where $\hat{t_i}$ means that the variable $t_i$ is omitted from the product.
In order to find its Grothendieck class, we first use the isomorphism between $X_\Gamma\setminus\Sigma_n$ and $X_{\Gamma^\vee}\setminus\Sigma_n$ to have
\begin{equation}
[X_\Gamma\setminus\Sigma_n]=[X_{\Gamma^\vee}\setminus\Sigma_n]\,.
\end{equation}
From this, we can find the class of $X_\Gamma\setminus\Sigma$ by studying $\mathcal{L}\cap\Sigma_n$\,. First, we need a result that tells us how to find the Grothendieck class of a hyperplane section of a given class in the Grothendieck ring.

\begin{lem}
Let $C$ be a class in the Grothendieck ring that can be written as a function of the torus class $\mathbb{T}$ by means of a polynomial expression $C=g(\mathbb{T})$. Then the transformation
\begin{equation}\label{hgtrans}
\mathcal{H}:g(\mathbb{T})\mapsto\frac{g(\mathbb{T})-g(-1)}{\mathbb{T}+1}
\end{equation}
gives an operation on the set of classes in the Grothendieck ring that are polynomial functions of the torus class $\mathbb{T}$ that can be interpreted as taking a hyperplane section.

\begin{proof}
Notice that for $\displaystyle g(\mathbb{T})=[\mathbb{P}^n]=\frac{(\mathbb{T}-1)^n-1}{\mathbb{T}}$, equation (\ref{hgtrans}) gives
\begin{equation}
\frac{g(\mathbb{T}-g(-1))}{\mathbb{T}+1}=\frac{(\mathbb{T}-1)^{n-1}-1}{\mathbb{T}}=[\mathbb{P}^{n-1}]\,,
\end{equation}
which effectively is the same as taking a hyperplane section. Since $\mathcal{H}$ is linear in $g$, it works for any $g$.
\end{proof}
\end{lem}

Since $[\Sigma_n]=\frac{(\mathbb{T}+1)^n-1-\mathbb{T}^n}{\mathbb{T}}=g(\mathbb{T})$, using the previous result gives that
\begin{equation}
[\mathcal{L}\cap\Sigma_n]=\frac{g(\mathbb{T})-g(-1)}{\mathbb{T}+1}=\frac{(1+\mathbb{T})^{n-1}-1}{\mathbb{T}}-\frac{\mathbb{T}^{n-1}-(-1)^{n-1}}{\mathbb{T}+1}\,,
\end{equation}
and with this, we find that
\begin{equation}
[X_{\Gamma^\vee}\setminus\Sigma_n]=[\mathcal{L}\setminus\Sigma_n]=[\mathcal{L}]-[\mathcal{L}\cap\Sigma_n]=\frac{\mathbb{T}^{n-1}-(-1)^{n-1}}{\mathbb{T}+1}\,.
\end{equation}

On the other hand, consider the variety $\mathcal{S}_n$ generated by the ideal
\begin{equation}
\mathcal{I}=(t_2t_3\cdots t_{n},t_1t_3\cdots t_n,\dots,t_1t_2\cdots\hat{t_i}\cdots t_n,\dots,t_1t_2\cdots t_{n-1})\,.
\end{equation}
This is the singularity subvariety of the divisor of singular normal crossings $\Sigma_n$ given by the union of coordinate hyperplanes (see \cite{a08} Section 1.3).

We can find the Grothendieck class of $\mathcal{S}_n$ by
\begin{lem}
The class of $\mathcal{S}_n$ is given by
\begin{equation}
[\mathcal{S}_n]=[\Sigma_n]-n\mathbb{T}^{n-2}=\sum_{i=2}^{n-1}\binom{n}{i}\mathbb{T}^{n-1-i}\,.
\end{equation}

\begin{proof}
Each coordinate hyperplane $\mathbb{P}_k^{n-2}$ in $\Sigma_n$ intersects the others along the union of its coordinate hyperplanes $\Sigma_{n-1}$. Thus, to obtain the class of $\mathcal{S}_n$ from the class of $\Sigma_n$, we just need to subtract the class of $n$ complements of $\Sigma_{n-1}$ in the $n$ components of $\Sigma_n$, from which the result follows.
\end{proof}
\end{lem}

Finally, this results lead to the Grothendieck class of the graph hypersurface of the banana graph by

\begin{thm}
The class of the graph hypersurface of the banana graph is given by
\begin{equation}
[X_{\Gamma^\vee}]=\frac{(\mathbb{T}+1)^n-1}{\mathbb{T}}-\frac{\mathbb{T}^n-(-1)^n}{\mathbb{T}}-n\mathbb{T}^{n-2}\,.
\end{equation}
\begin{proof}
Writing $[X_{\Gamma^\vee}]$ as
\begin{equation}
[X_{\Gamma^\vee}]=[X_{\Gamma^\vee}\setminus\Sigma]+[X_{\Gamma^\vee}\cap\Sigma]=[X_{\Gamma^\vee}\setminus\Sigma]+[\mathcal{S}_n]\,,
\end{equation}
and using the previous results for $[X_{\Gamma^\vee}\setminus\Sigma]$ and $[\mathcal{S}_n]$ we find the desired result.
\end{proof}
\end{thm}

\section{Irreducible Graph Hypersurfaces}

To start, notice that if a graph $\Gamma$ produces a reducible graph hypersurface $X_\Gamma$ so will its dual graph $\Gamma^\vee$.

\begin{thm}
 $X_\Gamma$ is irreducible if and only if $X_{\Gamma^\vee}$ is irreducible.

 \begin{proof}
This result follows directly from Proposition \ref{dualpol}.
 \end{proof}
\end{thm}

Even though every planar connected graph gives rise to a graph hypersurface, not every graph produces an irreducible graph hypersurface. To see this, consider a planar connected graph $\Gamma$ and following \cite{b08} define a \emph{separation} to be a decomposition of the graph into two subgraphs which share only one vertex. This common vertex is called a \emph{separating vertex}.

With this, we can give a relation between separating vertices and irreducibility of the graph hypersufaces by the following result

\begin{thm}\label{lemmared}
If the graph $\Gamma$ has a separating vertex such that each of the two components contain a cycle, then the corresponding graph hypersurface $X_\Gamma$ is reducible.

\begin{proof}
Let $a$ be a separating vertex in $\Gamma$ with the properties required. Consider the two associated components of the graph separately. Since the components do not share edges, the choice of edges to remove in order to achieve a spanning tree on each component is independent of one another. Thus every spanning tree of $\Gamma$ can be obtained by spanning trees of the two components joint through $a$. Then the graph polynomial for $\Gamma$ is the product of the graph polynomials of each component. Since each component contains a cycle, the corresponding graph polynomials are non-constant, and hence $X_\Gamma$ is reducible.
\end{proof}
\end{thm}

The concept of independence in the graph is the key for reducibility. With this discussion, we can refine the previous result by means of defining separability on graphs. A connected graph that has no separating vertices is called \emph{non-separable}, otherwise it is called \emph{separable}. The next result, due to Whitney, will provide the key ingredient to characterize reducibility on graph hypersurfaces by means of graph theoretical properties of $\Gamma$.

\begin{lem}\label{lemm52}
A connected graph is non-separable if and only if any two of its edges lie on a common cycle.

\begin{proof}
The proof can be found in Theorem 5.2 in \cite{b08}. The result follows from the definition of separating vertex and properties of cycles.
\end{proof}
\end{lem}

With this result, we can characterize reducibility of graph hypersurfaces.

\begin{thm}
$X_\Gamma$ is reducible if and only if $\Gamma$ is separable and each component contains a cycle.
\begin{proof}
By Theorem \ref{lemmared} we have that if $\Gamma$ is separable, and hence has a separating vertex, the corresponding graph polynomial $\Psi_\Gamma(t)$ is reducible.

Now, suppose that $\Gamma$ is such that $\Psi_\Gamma(t)$ is reducible. Thus
\begin{equation}
\Psi_\Gamma(t)=p(t)q(t)
\end{equation}
for $p(t)$ and $q(t)$ non constant polynomials. Now, let $t_i=x$ be variable and fix all other $t_j$, $j\neq i$. Therefore $\Psi_\Gamma(x)$ is a linear polynomial, and hence either $p(x)$ or $q(x)$ is a constant for any choice of $t_j$. Without loss of generality, suppose that $p(x)$ is a constant and that $q(x)$ is a linear polynomial. Now regard $p(x)$ and $q(x)$ as polynomials in $x$ whose coefficients depend on $t_j$. Since the coefficient of $x$ in $q(x)$ is non-zero for our initial choice of $t_j$, by continuity on $t_j$ it remains non-zero in a neighborhood of these initial values. Hence the coefficient of $x$ in $p(x)$ is zero for a neighborhood of $t_j$, and since $p(t)$ is a polynomial, $p(x)$ must be constant for all $t_j$.

Therefore, we have that $p(t)$ and $q(t)$ separate the variables $t_1,t_2,\dots, t_n$, that is, if $t_i$ is present in $p(t)$, then it is not in $q(t)$. Hence we have that $p(t)$ and $q(t)$ are homogeneous polynomials where the exponent of each $t_i$ is at most 1.

Since $p(t)$ and $q(t)$ separate the variables $t_i$ associated to the edges of $\Gamma$, the choice of the edges appearing in $p(t)$ is independent of the choice of the edges appearing in $q(t)$ in order to break cycles and obtain a spanning tree of $\Gamma$. Hence the edges appearing in $p(t)$ do not lie in any cycle in which are any of the edges that appear in $q(t)$.

Let
\begin{equation}
P=\{t_i| \text{ such that }t_i\text{ appears in }p(t)\}
\end{equation}
and
\begin{equation}
Q=\{t_i| \text{ such that }t_i\text{ appears in }q(t)\}\,.
\end{equation}

Since $P\cap Q=\emptyset$, $P\cup Q$ are all the variables appearing in $\Psi_\Gamma$, and no cycle in $\Gamma$ contains both elements of $P$ and $Q$, we have that, by Lemma \ref{lemm52}, $\Gamma$ is separable.
\end{proof}
\end{thm}

\end{document}